\documentclass [10pt,a4paper,draft] {article}
\usepackage [applemac] {inputenc}
\usepackage[french] {babel}
\usepackage {amsfonts}
\usepackage {amsmath}
\usepackage {amssymb}

\begin{document}

\begin{center}
\begin{LARGE}
Notes sur l'indice des alg\`ebres de Lie (I)
\end{LARGE}
\bigskip

par : Mustapha RA\"IS 
\footnotemark \footnotetext[1]{Universit\'e de POITIERS -
D\'epartement de Math\'ematiques - T\'el\'eport 2, Boulevard Marie et Pierre Curie - BP
30179 -
86962 FUTUROSCOPE CHASSENEUIL Cedex} 
\end{center}

\bigskip
\vspace{10 mm}

On trouvera ci-dessous : 

- Des d\'emonstrations de certains des r\'esultats obtenus r\'ecemment par D. Panyushev.

- Un compl\'ement portant sur ``l'in\'egalit\'e de Panyushev''.

- Des calculs d'indices de certaines contract\'ees d'alg\`ebres de Lie.

- Des exemples d'additivit\'e de l'indice des alg\`ebres de Lie.

- Des questions qui peuvent int\'eresser le lecteur.

\vspace{10 mm}
Je remercie Dmitri Panyushev pour l'impulsion donn\'ee par ses multiples travaux sur
l'indice des
alg\`ebres de Lie et pour la gentillesse avec laquelle il a r\'epondu \`a mes multiples
questions.

\bigskip
\underline{Notations} : Les notations sont en g\'en\'eral celles de \cite{Dix}. On se
reportera \`a
\cite{Ra} pour la d\'efinition de l'indice d'une repr\'esentation lin\'eaire d'une alg\`ebre
de Lie ; en
particulier :

\medskip
- Lorsque $\mathfrak{a}$ est un id\'eal d'une alg\`ebre de Lie $\mathfrak{g}$, l'indice
de la
repr\'esentation naturelle de $\mathfrak{g}$ dans $\mathfrak{a}$, sera not\'e
$\hbox{ind}(\mathfrak{g},
\mathfrak{a})$ comme dans \cite{Pa}. 

\medskip
- Lorsque $\mathfrak{h}$ est une sous-alg\`ebre de $\mathfrak{g}$, elle op\`ere
naturellement dans
l'espace vectoriel $\mathfrak{g}/\mathfrak{h}$, et on notera
$\hbox{ind}(\mathfrak{h},\mathfrak{g}/\mathfrak{h})$ l'indice de cette
repr\'esentation de
$\mathfrak{h}$.

\vspace{15mm}

%%%%%%%%%%%%%%%%%%%%%%%
\section{L'in\'egalit\'e de Panyushev}
        
\noindent
\textbf{1.1.}~On aura besoin du r\'esultat \'el\'ementaire et classique suivant : soient
$\mathfrak{g}$
une alg\`ebre de Lie, $\mathfrak{a}$ un id\'eal de $\mathfrak{g}$, $\ell$ une forme
lin\'eaire sur
$\mathfrak{g}$, $\ell_0$ la restriction de $\ell$ \`a $\mathfrak{a}, \mathfrak{h} =
\mathfrak{a}^\ell$ l'orthogonal de $\mathfrak{a}$ relativement \`a la forme $B_\ell$
associ\'ee \`a
$\ell, \ell_1$ la restriction de $\ell$ \`a $\mathfrak{h}$. Alors :
$\mathfrak{h}^{\ell_1} =
\mathfrak{a}^{\ell_0}+ \mathfrak{g}^\ell$. (En effet : $\mathfrak{h}^\ell =
(\mathfrak{a}^\ell)^\ell = \mathfrak{a} + \mathfrak{g}^\ell$ et
$\mathfrak{h}^{\ell_1} =
\mathfrak{h} \cap \mathfrak{h}^\ell = \mathfrak{a} \cap \mathfrak{h} +
\mathfrak{g}^\ell =
\mathfrak{a}^{\ell_0} + \mathfrak{g}^\ell)$. (Voir \cite{Dix}, 1.12.4.(i).)

\vskip 7mm
\noindent
\textbf{1.2.}~Dans \cite{Pa}, Panyushev d\'emontre :

        1.4. \textbf{Th\'eor\`eme} : \textit{Soit $\mathfrak{a}$ un id\'eal dans une alg\`ebre de Lie
$\mathfrak{g}$. On a alors (ce qu'on appellera dans la suite \underline{l'in\'egalit\'e
de Panyushev}) :
}$\hbox{ind}(\mathfrak{g}) + \hbox{ind}(\mathfrak{a}) \leq
\dim(\mathfrak{g}/\mathfrak{a}) + 2
\, \hbox{ind}(\mathfrak{g},\mathfrak{a})$.

\medskip
\noindent
\textsc{Une d\'emonstration} : Soient $\ell$ dans $\mathfrak{g}^*, \ell_0 =
\ell_{|\mathfrak{a}},\
\mathfrak{h} = \mathfrak{a}^\ell$, de sorte que : $\mathfrak{a}^{\ell_0} +
\mathfrak{g}^\ell
\subset \mathfrak{h}$ (pr\'ecis\'ement $\mathfrak{a}^{\ell_0} + \mathfrak{g}^\ell =
\mathfrak{h}^{\ell_1}$, avec $\ell_1 = \ell|_\mathfrak{h}$). On \'ecrit simplement :
$\dim(a^{\ell_0} + \mathfrak{g}^\ell) \leq \dim \, \mathfrak{h}$. 

\bigskip
        . Sachant que $\mathfrak{a}^{\ell_0} \cap \mathfrak{g}^\ell = \mathfrak{a} \cap
\mathfrak{g}^\ell$ est l'orthogonal (dans $\mathfrak{a}$) du sous-espace
$\mathfrak{g}.\ell_0$
de $\mathfrak{a}^*$, il vient $\dim(\mathfrak{a}^{\ell_0}\cap \mathfrak{g}^\ell) =
\hbox{codim}(\mathfrak{g}.\ell_0)$.

\bigskip
        . On a : $\dim\, \mathfrak{h} = \dim\, \mathfrak{a}^\ell = \dim\, \mathfrak{g} -
\dim\,
\mathfrak{a} + \dim(\mathfrak{a} \cap \mathfrak{g}^\ell) =
\dim(\mathfrak{g}/\mathfrak{a}) +
\hbox{codim}(\mathfrak{g}.\ell_0)$.

\medskip
On a donc : $\dim\, \mathfrak{a}^{\ell_0} + \dim\, \mathfrak{g}^\ell \leq
\dim(\mathfrak{g}/\mathfrak{a}) + 2\, \hbox{codim}(\mathfrak{g}.\ell_0)$.

        D'o\`u : $\hbox{ind}(\mathfrak{g}) + \hbox{ind}(\mathfrak{a}) \leq
\dim(\mathfrak{g}/\mathfrak{a}) + 2\,
\hbox{codim}(\mathfrak{g}.\ell_0)$ et ce, pour toute forme lin\'eaire $\ell_0$ sur
$\mathfrak{a}$.
Par suite :
$$
        \hbox{ind}(\mathfrak{g}) + \hbox{ind}(\mathfrak{a}) \leq
\dim(\mathfrak{g}/\mathfrak{a}) + 2\,
\hbox{ind}(\mathfrak{g}, \mathfrak{a}).
$$

\vskip 12mm

%%%%%%%%%%%%%%%%%%%%%%%
\section{L'indice du normalisateur d'un centralisateur}

\noindent
\textbf{2.1.}~Soit $e$ un \'el\'ement nilpotent d'une alg\`ebre de Lie simple
$\mathfrak{g}$ (le corps
de base est, pour moi, le corps des complexes). On note $z$ le centralisateur de $e$
(dans
$\mathfrak{g})$, $\eta$ le normalisateur de $z$ (dans $\mathfrak{g}$) et $\delta$ le
centre de
$z$.  Dans \cite{Pa}, Panyushev d\'emontre :

\medskip
        4.4. (ii)~\textbf{Th\'eor\`eme}~:~\textit{Sous l'hypoth\`ese} ``$\hbox{ind}(\eta,
\delta)=0$'', \textit{on
a : }$\hbox{ind}(\eta) = \hbox{ind}(z) - \dim\, \delta$ \textit{\underline{et}}
$\hbox{ind}(\eta, z)
= \hbox{ind}(\eta)$.

\medskip
\noindent
\textsc{Une d\'emonstration} : Soit $\ell$ une forme lin\'eaire sur $\eta$. On applique
1.1. avec
$\mathfrak{g} = \eta \ \  \mathfrak{a} = \delta$ ; avec $\mathfrak{h} = \delta^\ell$ et
$\ell_1 = \ell_{|\mathfrak{h}}$, on a donc : $\mathfrak{h}^{\ell_1} = \delta +
\eta^\ell$
($\delta^{\ell_0}= \delta$ puisque $\delta$ est ab\'elien) et clairement :
$\mathfrak{h}^{\ell_1} =
\mathfrak{h}^\ell$. 

\bigskip
- Il est imm\'ediat que : ``$\delta \cap \eta^\ell = \{0\}$'' ssi  ``$\eta.\ell_0 =
\delta^*$'', i.e. ssi
``$\hbox{ind}(\eta,\delta)=0$ \underline{et} $\ell_0$ est un \'el\'ement r\'egulier de
$\delta^*$''.

\bigskip
- Notons $B : \delta \times \eta \longrightarrow \mathbb{C}$ la forme d\'efinie par : 
$B(x,y) =\ 
<\ell, [x,y]>$. On voit que $B(x,y) = 0$ pour tout $y$ dans $\eta$ \'equivaut \`a : $x
\in \delta \cap
\eta^\ell$, et $B(x,y) = 0$ pour tout $x$ dans $\delta$ \'equivaut \`a : $y \in
\delta^\ell$. Comme
$\delta$ et $\eta/z$ ont m\`eme dimension (\cite{B-K}), il vient que les assertions
suivantes sont
\'equivalentes :
$$
\eta.\ell_0 = \delta^* ; B \ \hbox{induit une dualit\'e entre} \ \delta \ \hbox{et}\ 
(\eta/z)\,  ;\ 
\delta^\ell = z.
$$

\bigskip
- Supposons donc $\eta.\ell_0 = \delta^*$. Avec $\ell_1 = \ell_{|z}$, on a alors : 
$$
        z^{\ell_1} = \delta \oplus \eta^\ell \quad \hbox{et}\quad  z^{\ell_1} = z^\ell.
$$
\indent
\quad . De $z^{\ell_1} = \delta \oplus \eta^\ell$, il r\'esulte : $\dim\ z^{\ell_1} =
\dim\ \delta +
\dim\, \eta^\ell$.

\medskip
\indent
\quad . De $z^\ell = z^{\ell_1}$, il r\'esulte : $\dim\ \eta.\ell_1 = \dim\ \eta -
\dim\ z^\ell = \dim\
z + \dim\ \delta - (\dim\ \delta + \dim\ \eta^\ell) = \dim\ z - \dim\ \eta^\ell$, d'o\`u
$\hbox{codim}(\eta.\ell_1) = \dim \ \eta^\ell$. En r\'esum\'e : 
$$
        \dim \ z^{\ell_1} = \dim \ \delta + \dim \ \eta^\ell
$$
$$
        \hbox{codim}(\eta.\ell_1) = \dim \ \eta^\ell.
$$
\indent
        Le corps de base \'etant le corps des complexes, on d\'eduit des deux \'egalit\'es
pr\'ec\'edentes :
$$
        \hbox{ind}\ z = \dim\ \delta + \hbox{ind}(\eta)
$$
$$
        \hbox{ind}(\eta,z) = \hbox{ind}\ \eta
$$

\vskip 7mm
\noindent
\textbf{2.2.}~\textsc{Remarques} :  $\centerdot$ La d\'emonstration de Panyushev
utilise l'in\'egalit\'e
du paragraphe 1.

\smallskip
\noindent
$\centerdot$ Toujours sous l'hypoth\`ese ``$\eta.\ell_0 = \delta^*$'', les relations
trouv\'ees
ci-dessus liant
$z^{\ell_1},\ \eta^\ell$ et $\eta.\ell_1$ montrent que les 3 assertions suivantes
sont 2 \`a 2
\'equivalentes :

\medskip
\indent
        . $\ell$ est un \'el\'ement r\'egulier de $\eta^*$ (i.e. $\dim\ \eta^\ell =
\hbox{ind}(\eta)$)

\smallskip
\indent
        . $\ell_1$ est un \'el\'ement r\'egulier de $z^*$ (i.e. $\dim\ z^{\ell_1} = \hbox{ind}\ z$)

\smallskip
\indent
. $\ell_1$ est un \'el\'ement r\'egulier de $z^*$, pour la repr\'esentation naturelle de
$\eta$ dans
$z^*$\break(i.e. $\hbox{codim}(\eta.\ell_1) = \hbox{ind}(\eta,z))$. 

\vskip 7mm
\noindent
\textbf{2.3.}~On se permettra d'exprimer la condition ``$\hbox{ind}(\eta,\delta) =
0$'' en disant que
$\eta$ admet une orbite ouverte dans $\delta^*$. Comme expliqu\'e ci-dessus, Panyushev
d\'emontre les formules des indices de $\eta$, de $z,\ldots ,$ sous cette condition.

\vskip 7mm
\noindent
\textsc{Question 1} : Est-ce que ``$\hbox{ind}(\eta,\delta) = 0$'' \'equivaut \`a
``$\hbox{ind}(\eta) 
= \hbox{ind}(z)- \dim\
\delta$ \underline{et} : $\hbox{ind}(\eta, z) = \hbox{ind}\ \eta$'' ?

\vskip 7mm
\noindent
\textbf{2.4.}~On trouve dans \cite{Pa} une liste de cas o\`u il est prouv\'e que la
condition
``$\hbox{ind}(\eta,\delta) = 0$'' est v\'erifi\'ee (y voir le th\'eor\`eme 4.7). Il faut
noter que dans tous 
ces cas, Panyushev d\'emontre une propri\'et\'e plus forte que celle de l'orbite ouverte,
\`a savoir : 
$\eta$ admet dans $\delta^*$ un \underline{nombre fini} d'orbites.

\vskip 7mm
\noindent
\textbf{2.5.}~Dans le cas o\`u $\mathfrak{g} = sl(n)$, on pr\'esente ci-dessous une
d\'emonstration
\underline{tr\`es peu diff\'erente} de celle de Panyushev. Soit $e$ un \'el\'ement nilpotent
non nul de
$sl(n)$ ; le centre $\delta$ du centralisateur $z$ de $e$ admet $(e,e^2,\ldots
,e^r)$ comme base,
$r$ \'etant l'entier v\'erifiant : $r \geq 1,\ e^r \not= 0,\ e^{r+1} = 0$. Soit $h$ dans
$sl(n)$ tel que
$[h,e]= 2e$ ; alors $[he^k, e] = 2e^{k+1}$ (ici $he^k$ est le produit effectu\'e dans
l'alg\`ebre
associative des matrices), et $[he^k,e^\ell] = 2\ell e^{k+\ell}$. Quitte \`a projeter
(s'il le faut) les
$he^k$ dans $sl(n)$, parall\`element \`a l'espace des matrices scalaires, on peut
affirmer : il existe
$x_0, x_1,\ldots ,x_{r-1}$ dans $sl(n)$ tels que : $[x_k, e^\ell] = 2\ell
e^{k+\ell}\ (0 \leq k \leq
r-1,\ 1 \leq \ell \leq r)$ ; par suite, le normalisateur $\eta$ de $z$ dans $sl(n)$
s'\'ecrit : $\eta =
z\oplus \mathfrak{q}$, avec : $\displaystyle \mathfrak{q} = \sum^{r-1}_{k=0}\
\mathbb{C}\
x_k$. Soit $\ell_0$ une forme lin\'eaire sur $\delta$ ; on pose $\lambda_k = \
<\ell_0, e^k>\ (1 \leq
k \leq r)$, de sorte que $\displaystyle \ell_0 = \sum^r_1\, \lambda_k\ \xi_k$ o\`u
$(\xi_1,\ldots
, \xi_r)$ est la base duale de $(e, e^2,\ldots , e^r)$. On cherche la matrice $A$ de
l'application
lin\'eaire : $x \mapsto x.\ell_0$, de $\mathfrak{q}$ dans $\delta^*$, relativement aux
bases
$(x_0,x_1,\ldots ,x_{r-1})$ de $\mathfrak{q}$ et $(\xi_1,\ldots ,\xi_r)$ de
$\delta^*$ : 
$$
        a_{ij} = \ <x_j.\ell_0, e^i> \ = -<\ell_0, [x_j,e^i]>\  = - <\ell_0, 2ie^{i+j}>
$$
d'o\`u : $a_{ij} = 0$ lorsque $i>r-j$ et $a_{ij} = -2i\ \lambda_{i+1}$ lorsque $i+j
\leq r$. Il vient
alors : $|\det A| = 2^r(r!)\lambda^r_r$.

\smallskip
        Donc $\hbox{ind}(\eta, \delta) = 0$ et $\ell_0$ est une forme r\'eguli\`ere, au sens :
$\eta.\ell_0 =
\delta^*$, si et seulement si : $<\ell_0,e^r> \ \not= 0$.

\smallskip
        (A quelques modifications pr\`es, cette d\'emonstation -(encore une fois, c'est celle
de Panyushev,
habill\'ee une peu diff\'eremment)- vaut pour les alg\`ebres de Lie $so(2n+1), sp_{2n}$ et
pour certains
types d'\'el\'ements nilpotents dans $so(2n)$, voir le th\'eor\`eme 4.7 de \cite{Pa}.)

\medskip
        $\centerdot$ On va voir que dans cet exemple, on peut calculer l'indice de
l'alg\`ebre $\eta/z$. On a :
$[he^k, he^\ell] = 2(\ell-k)h\ e^{k+\ell}$ ($=0$ lorsque $k+\ell \geq r+1$). De plus
$[h\ e^r,e] =
2\ e^{r+1} = 0$, d'o\`u $h\ e^r \in z$. De ceci il r\'esulte que $(\bar{x}_0,
\bar{x}_1,\ldots , \bar{x}_{r-1}$, (o\`u $\bar{x}_j = x_j \mod z)$ est une base de
$\eta/z$, avec les
crochets :

\vfill\eject
%\bigskip
\quad $[\bar{x}_k, \bar{x}_\ell] = 2(\ell-k)\bar{x}_{k+\ell}$ lorsque $k+\ell \leq r-1$

\smallskip
        \qquad \qquad \ $ = 0$ dans les autres cas.

%$$
%\begin{eqnarray}
%        [\bar{x}_k, \bar{x}_\ell]
%        &= &2(\ell-k)\bar{x}_{k+\ell}\ \mathrm{lorsque}\ k+\ell \leq r-1\nonumber\\
%        &=        &0 \ \mathrm{dans les autres cas}
%\end{eqnarray}
%$$

\bigskip
        Soit $\lambda$ une forme lin\'eaire sur $\eta/z,\ \lambda = \displaystyle \sum^{r-1}_0\
\lambda_k\ \eta_k$, o\`u $(\eta_0,\ldots , \eta_{r-1}$) est la base duale de
$(\bar{x}_0,\ldots
,\bar{x}_{r-1})$, et soit $A$ la matrice de l'application $x \mapsto x.\lambda$ de
$\eta/z$ dans
son dual, relativement aux bases $(\bar{x}_0,\ldots ,\bar{x}_{r-1})$ et
$(\eta_0,\ldots ,\eta_r)$ : 

\bigskip
\quad $a_{ij} = \ <\bar{x}_j.\lambda,\ \bar{x}_i>\  = \ <\lambda, [\bar{x}_i,\
\bar{x}_j]> \ =0$
lorsque $i+j \geq r$

\smallskip
\quad $a_{ij} = 2(j-i)\lambda_{i+j}$ lorsque $i+j \leq r-1$.

\bigskip
A nouveau : $\det A = c\, \lambda^r_{r-1}$ avec $c\in \mathbb{C}$, et $c \not= 0$ si et
seulement si $r$ est \underline{pair}. Lorsque $r$ est impair, il y a un et un seul
\'el\'ement sur
l'anti-diagonale de $A$ qui soit nul, et il vient que $A$ admet un mineur non nul de
taille $(r-1)$.

\bigskip
\underline{Conclusion} : L'alg\`ebre de Lie $\eta/z$ est d'indice z\'ero (resp. d'indice
1) lorsque sa
dimension est paire (resp. impaire).

\vskip 7mm
\noindent
\textsc{Question 2} : Quel est l'indice de $\eta/z$ (au moins lorsque
$\hbox{ind}(\eta,\delta) = 0$) ?

\bigskip
$\centerdot$ Dans les exemples matriciels trait\'es ci-dessus, on constate que, aussi
bien pour la
repr\'esenta\-tion naturelle de $\eta$ dans $\delta$ que pour la repr\'esentation
adjointe de
$\eta/z$, les formes lin\'eaires r\'eguli\`eres sont caract\'eris\'ees par le fait
qu'\underline{elles ne
s'annulent pas sur l'\'el\'ement de plus grande} $h$\underline{-gra\-dua\-tion}. On
dispose par
ailleurs d'une application lin\'eaire $ad(e) : \eta \longrightarrow \delta$, dont le
noyau est $z$, de
sorte que cette application induit un isomorphisme d'espaces vectoriels de $\eta/z$
sur $\delta$,
qu'on continuera \`a noter $ad(e)$. On voit alors que, lorsque $\ell_0$ est une forme
lin\'eaire sur
$\delta$, on a l'\'equivalence des 2 assertions suivantes :

\medskip
        - $\ell_0$ est un \'el\'ement r\'egulier au sens de : $\eta.\ell_0 = \delta^*$

\medskip
        - $\ell_0 \, o \,  ad(e)$ est un \'el\'ement r\'egulier de $(\eta/z)^*$.

\vskip 12mm

%%%%%%%%%%%%%%%%%%%%%%%
\section{Des cas d'\'egalit\'e dans l'in\'egalit\'e de Panyushev}
        
\noindent
\textbf{3.1.}~On reprend les notations du paragraphe 1 : $\mathfrak{g}$ est une
alg\`ebre de Lie,
$\mathfrak{a}$ en est un id\'eal, $\ell$ est une forme lin\'eaire sur $\mathfrak{g}$ et
$\ell_0$ est sa
restriction \`a $\mathfrak{a}$. On a alors, de fa\`eon \'evidente : $\mathfrak{a}^{\ell_0} +
\mathfrak{g}^\ell \subset \mathfrak{a}^\ell$ et l'in\'egalit\'e qui exprime que la
dimension de
l'espace vectoriel $\mathfrak{a}^{\ell_0} + \mathfrak{g}^\ell$ est major\'ee par celle de
$\mathfrak{a}^\ell$ s'\'ecrit : 
$$
        \dim \mathfrak{a}^{\ell_0} + \dim \mathfrak{g}^\ell \leq
\dim(\mathfrak{g}/\mathfrak{a}) + 2
\ \hbox{codim}(\mathfrak{g}.\ell_0).
$$
On en a d\'eduit l'in\'egalit\'e de Panyushev. On s'int\'eresse aux \underline{cas
d'\'egalit\'e} dans cette
in\'egalit\'e.

\vskip 7mm
\noindent
\textbf{Lemme} :  \textit{Notons $\Omega$ l'ouvert de Zariski dans $\mathfrak{g}^*$
constitu\'ee
par les formes lin\'eaires $\ell$ telles que} $\hbox{codim}(\mathfrak{g}.\ell_0) = \,
\hbox{ind}(\mathfrak{g},\mathfrak{a})$ \textit{(i.e. telles que $\ell_0 =
\ell_{|\mathfrak{a}}$ soit
un \'el\'ement r\'egulier pour la repr\'esentation naturelle de $\mathfrak{g}$ dans
$\mathfrak{a}$).}

        \textit{L'\'egalit\'e :} $\hbox{ind} \, \mathfrak{a} + \hbox{ind} \, \mathfrak{g} =
\dim(\mathfrak{g}/\mathfrak{a}) + 2 \ \hbox{ind}\, (\mathfrak{g},\mathfrak{a})$
\textit{est
r\'ealis\'ee si et seulement si :}
$$
        \hbox{\textit{Pour toute}} \ \ell \ \hbox{\textit{dans}} \ \Omega, \hbox{\textit{on
a :}} \
\mathfrak{a}^{\ell_0} + \mathfrak{g}^\ell = \mathfrak{a}^\ell.
$$

\noindent
\textsc{D\'emonstration} : $\centerdot$ Supposons r\'ealis\'ee l'\'egalit\'e de Panyushev.
Soit $\ell$ dans
$\Omega$. On a alors : $\dim \mathfrak{a}^{\ell_0} + \dim \mathfrak{g}^\ell \leq
\dim(\mathfrak{g}/\mathfrak{a}) + 2\ \hbox{ind}(\mathfrak{g},\mathfrak{a})$ et
$\hbox{ind}\
\mathfrak{a} +\ \hbox{ind}\ \mathfrak{g} = \dim(\mathfrak{g}/\mathfrak{a}) + 2\
\hbox{ind}(\mathfrak{g},\mathfrak{a}) \leq \dim \mathfrak{a}^{\ell_0} + \dim
\mathfrak{g}^\ell \leq \dim(\mathfrak{g}/\mathfrak{a}) + 2\
\hbox{ind}(\mathfrak{g},\mathfrak{a})$.

\bigskip
        Donc $\dim \mathfrak{a}^{\ell_0} + \dim \mathfrak{g}^\ell = \dim
\mathfrak{g}/\mathfrak{a}
+ 2\ \hbox{ind}(\mathfrak{g}/\mathfrak{a})$ (d'o\`u r\'esulte d'ailleurs que $\ell$ et
$\ell_0$ sont
des formes lin\'eaires r\'eguli\`eres respectivement sur les alg\`ebres de Lie
$\mathfrak{a}$ et
$\mathfrak{g}$), et ainsi :
$$
        \dim \mathfrak{a}^{\ell_0} + \dim \mathfrak{g}^\ell =
\dim(\mathfrak{g}/\mathfrak{a}) + 2\
\hbox{codim}(\mathfrak{g}.\ell_0)
$$
c'est-\`a-dire : $\dim(\mathfrak{a}^{\ell_0} + \mathfrak{g}^\ell) = \dim
\mathfrak{a}^\ell$.

\medskip
\noindent
$\centerdot$ Supposons qu'on ait, pour tout $\ell$ dans $\Omega$ : 
$$
        \mathfrak{a}^{\ell_0} + \mathfrak{g}^\ell = \mathfrak{a}^\ell
$$
c'est-\`a-dire : $\dim \mathfrak{a}^{\ell_0} + \dim \mathfrak{g}^\ell = \dim
\mathfrak{g}/\mathfrak{a} + 2\ \hbox{ind}(\mathfrak{g}, \mathfrak{a})$.

\smallskip
        Dans l'ouvert de Zariski $\Omega$, on peut trouver $\ell$ qui soit une forme
r\'eguli\`ere sur
$\mathfrak{g}$ et telle que $\ell_0 = \ell_{|\mathfrak{a}}$ soit une forme r\'eguli\`ere
sur
$\mathfrak{a}$. Donc : $\hbox{ind}\ \mathfrak{a} + \hbox{ind}\ \mathfrak{g} = \dim 
\mathfrak{g}/\mathfrak{a} + 2\ \hbox{ind}(\mathfrak{g}, \mathfrak{a})$.

\vskip 7mm
 %%%%
\noindent
\textbf{3.2.}~On se restreint, jusqu'\`a mention explicite du contraire, au cas o\`u
l'id\'eal
$\mathfrak{a}$ est \underline{ab\'elien}, o\`u il est plus facile de reconna\^\i tre les cas
d'\'egalit\'e.

\bigskip
        On notera dans ce cas que $\mathfrak{a}^{\ell_0} = \mathfrak{a}$, et que la
sous-alg\`ebre
$\mathfrak{a} + \mathfrak{g}^\ell$ est \underline{subordonn\'ee} \`a $\ell$, tandis que
$\mathfrak{a}^\ell$ est \underline{co-isotrope}. 

\vskip  7mm
\noindent
\textbf{Lemme}  : \textit{Les assertions suivantes sont \'equivalentes :}

\smallskip
\noindent
(1) $\hbox{ind}(\mathfrak{g}) + \dim \mathfrak{a} = \dim(\mathfrak{g}/\mathfrak{a}) +
2\ \hbox{ind}(\mathfrak{g},\mathfrak{a})$.

\smallskip
\noindent
(2) \begin{itshape} Pour toute $\ell$ dans $\Omega$, on a : $\mathfrak{a} +
\mathfrak{g}^\ell =
\mathfrak{a}^\ell$.\end{itshape}

\smallskip
\noindent
(3) \begin{itshape}Pour toute $\ell$ dans $\Omega$, $\mathfrak{a} +
\mathfrak{g}^\ell$ est une
\underline{polarisation} en $\ell$.\end{itshape}

\smallskip
\noindent
(4) \begin{itshape}Pour toute $\ell$ dans $\Omega$,  $\mathfrak{a}^\ell$ est une
\underline{polarisation} en
 $\ell$.\end{itshape}

\vskip 5mm
\noindent
\textsc{D\'emonstration} : $\centerdot$ L'\'equivalence de (1) et (2) a \'et\'e prouv\'ee.

\medskip
\noindent
$\centerdot$ Soit $\ell$ dans $\mathfrak{g}^*$ telle que $\mathfrak{a} +
\mathfrak{g}^\ell =
\mathfrak{a}^\ell$. Alors $(\mathfrak{a} + \mathfrak{g}^\ell)^\ell =
(\mathfrak{a}^\ell)^\ell =
\mathfrak{a} + \mathfrak{g}^\ell$. Donc $\mathfrak{a} + \mathfrak{g}^\ell$ est une
polarisation
en $\ell$.

\medskip
\noindent
$\centerdot$ Soit $\ell$ dans $\mathfrak{g}^*$ telle que $(\mathfrak{a} +
\mathfrak{g}^\ell)$
soit une polarisation en $\ell$, i.e. telle que $\mathfrak{a} + \mathfrak{g}^\ell =
(\mathfrak{a} +
\mathfrak{g}^\ell)^\ell$ ; comme $(\mathfrak{a} + \mathfrak{g}^\ell)^\ell =
\mathfrak{a}^\ell$,
on a : $\mathfrak{a} + \mathfrak{g}^\ell = \mathfrak{a}^\ell$, ce qui s'\'ecrit :
$\mathfrak{a}^\ell
= (\mathfrak{a}^\ell)^\ell$.

\medskip
\noindent
$\centerdot$ Soit $\ell$ dans $\mathfrak{g}^*$ telle que $\mathfrak{a}^\ell$ soit une
polarisation en $\ell$. Alors $\mathfrak{a}^\ell = (\mathfrak{a}^\ell)^\ell =
\mathfrak{a} +
\mathfrak{g}^\ell$.

\bigskip
        Ceci montre que les conditions (2), (3), (4) sont 2 \`a 2 \'equivalentes.

\vskip 9mm
 %%%%
\noindent
\textbf{3.3.}~On examine ici plus particuli\`erement le cas o\`u $\mathfrak{a}$ est un
id\'eal ab\'elien
qui est \underline{facteur direct}, i.e. on suppose qu'il existe une sous-alg\`ebre
$\mathfrak{q}$ de
$\mathfrak{g}$ telle que $\mathfrak{g} = \mathfrak{q} \oplus \mathfrak{a}$.
Lorsqu'il en est
ainsi, on a : $\mathfrak{a}^\ell = \mathfrak{a} \oplus \mathfrak{q}(\ell_0)$, o\`u
$\mathfrak{q}(\ell_0)$ est l'annulateur de $\ell_0$ dans $\mathfrak{q}$, i.e.
l'ensemble des $x$
dans $\mathfrak{q}$ tels que : $<\ell, [x,\mathfrak{a}]> \ = 0$, et ce pour toute
forme lin\'eaire
$\ell$.

\vfill\eject
%\vskip 7mm
\noindent
\textbf{Lemme}  : \textit{Les assertions suivantes sont \'equivalentes :}

\smallskip
\noindent
(1) $\hbox{ind}\ \mathfrak{g} + \dim \mathfrak{a} = \dim(\mathfrak{g}/\mathfrak{a})
+ 2\
\hbox{ind}(\mathfrak{g},\mathfrak{a})$.

\smallskip
\noindent
(2) \begin{itshape} Pour toute $\ell$ dans $\Omega$, le ``stabilisateur g\'en\'erique''
$\mathfrak{q}(\ell_0)$ est ab\'elien. \end{itshape}

\vskip 5mm
\noindent
\textsc{D\'emonstration} : $\centerdot$ On suppose que (1) est v\'erifi\'ee.  Soit $\ell$
dans $\Omega$ ;
alors
$\mathfrak{a}^\ell =
\mathfrak{a} + \mathfrak{g}^\ell$, i.e. : $\mathfrak{a}^\ell = \mathfrak{a} +
\mathfrak{q}(\ell_0)$ ; comme
$\ell$ est r\'eguli\`ere, $\mathfrak{g}^\ell$ est ab\'elienne, d'o\`u : $[\mathfrak{a}^\ell,
\mathfrak{a}^\ell] \subset [\mathfrak{a}, \mathfrak{g}^\ell] \subset \mathfrak{a}$ ;
en utilisant
l'\'egalit\'e $\mathfrak{a}^\ell = \mathfrak{a} + \mathfrak{q}(\ell_0)$, il vient :
$[\mathfrak{q}(\ell_0),\mathfrak{q}(\ell_0)] \subset \mathfrak{a}$, d'o\`u
$\mathfrak{q}(\ell_0)$
est ab\'elienne.

\medskip
\noindent
$\centerdot$ On suppose que (2) est v\'erifi\'ee. Soit $\ell_0$ dans $\mathfrak{a}^*$
telle que
$\mathfrak{q}(\ell_0)$ soit ab\'elienne, et soit $\ell$ dans $\mathfrak{g}^*$, telle
que $\ell_0 =
\ell_{|\mathfrak{a}}$. On a alors : 
$$
        <\ell, [\mathfrak{a}+ \mathfrak{q}(\ell_0), \mathfrak{a} + \mathfrak{q}(\ell_0)]> \
= \ <\ell,
[\mathfrak{a}, \mathfrak{q}(\ell_0)]>\  = 0
$$
par d\'efinition de $\mathfrak{q}(\ell_0)$. Ceci exprime que $\mathfrak{a} +
\mathfrak{q}(\ell_0)
= \mathfrak{a}^\ell$ est une sous-alg\`ebre subordonn\'ee \`a la forme lin\'eaire $\ell$.
Comme par
ailleurs $\mathfrak{a}^\ell$ est co-isotrope relativement \`a $\ell$
$((\mathfrak{a}^\ell)^\ell =
\mathfrak{a }+ \mathfrak{g}^\ell \subset \mathfrak{a}^\ell)$, on a n\'ecessairement
$\mathfrak{a}^\ell = (\mathfrak{a}^\ell)^\ell = \mathfrak{a} + \mathfrak{g}^\ell$,
i.e. :
$\mathfrak{a}^\ell$ est une polarisation en $\ell$.

\medskip
\noindent
$\centerdot$ Ce lemme permet de donner un certain nombre d'exemples o\`u il y a \'egalit\'e :

\medskip
        \begin{enumerate}
        \item Soit $\mathfrak{q}$ une alg\`ebre de Lie simple et soit $\rho : \mathfrak{q}
\longrightarrow
\mathfrak{g}\ell(V)$ une repr\'esentation lin\'eaire de dimension finie. Lorsque cette
dimension est
strictement sup\'erieure \`a celle de $\mathfrak{q}$, il est connu que le stabilisateur
g\'en\'erique est
r\'eduit \`a $\{0\}$. Il y a donc \'egalit\'e dans l'in\'egalit\'e de Panyushev pour
$\mathfrak{g} = V
\underset{\rho}{\ltimes} \mathfrak{q}$, avec $\mathfrak{a} = V \times \{0\}$. 

        \medskip
        \item Soit $\mathfrak{q}$ une alg\`ebre de Lie et soit : $\rho : \mathfrak{q}
\longrightarrow
\mathfrak{g}\ell(\mathfrak{q})$ sa repr\'esentation adjointe. On d\'esigne par
$\mathfrak{g} =
\mathfrak{q} \underset{\rho}{\ltimes} \mathfrak{q}$ l'alg\`ebre de Lie produit
semi-direct associ\'e \`a
cette repr\'esentation. On sait (\cite{Ra}) que : $\hbox{ind}\ \mathfrak{g} = 2\
\hbox{ind}\
\mathfrak{q}$. On a donc, avec $\mathfrak{a} = \mathfrak{q} \times \{0\}$ : 
$$
        \hbox{ind}\ \mathfrak{g} + \dim \mathfrak{a} = 2\ \hbox{ind}\ \mathfrak{q} + \dim
\mathfrak{q}
$$
$$
        \hbox{et}\ \dim(\mathfrak{g}/\mathfrak{a}) + 2\
\hbox{ind}(\mathfrak{g},\mathfrak{a}) = \dim
\mathfrak{q} + 2\ \hbox{ind}\ \mathfrak{q}.
$$

Il y a donc \'egalit\'e, ce qui peut se voir a priori en notant que le stabilisateur
g\'en\'erique de la
repr\'esentation coadjointe d'une alg\`ebre de Lie est ab\'elien.
        \end{enumerate}

\vskip 7mm
 %%%%
\noindent
\textbf{3.4.}~On notera enfin une autre caract\'erisation du cas d'\'egalit\'e, toujours
lorsque
$\mathfrak{a}$ est ab\'elien, qui est susceptible d'avoir une interpr\'etation g\'eom\'etrique.

\vskip 7mm
\noindent
\textbf{Lemme}  : \textit{Les assertions suivantes sont \'equivalentes :}

\smallskip
\noindent
(1) $\hbox{ind}\ \mathfrak{g} + \dim \mathfrak{a} = \dim(\mathfrak{g}/\mathfrak{a})
+ 2\
\hbox{ind}(\mathfrak{g},\mathfrak{a})$.

\smallskip
\noindent
(2) \begin{itshape} Pour toute $\ell$ dans $\Omega$, on a : $\dim
(\mathfrak{g}.\ell) = 2\
\dim(\mathfrak{g}.\ell_0)$. \end{itshape}

\bigskip
En effet, l'\'egalit\'e (1) s'\'ecrit : 
$$
        \dim \mathfrak{g} - \hbox{ind} \ \mathfrak{g} =2(\dim \mathfrak{a} -
\hbox{ind}(\mathfrak{g},\mathfrak{a})).
$$

\vskip 7mm
 %%%%
\noindent
\textbf{3.5.}~Il y a des exemples int\'eressants o\`u $\mathfrak{a}$ \underline{n'est
pas ab\'elien} et
o\`u il y a \'egalit\'e dans l'in\'egalit\'e de Panyushev. On en donne ci-dessous une liste
probablement
incompl\`ete.

        \begin{enumerate} 
        \item Soit $\mathfrak{g}$ une alg\`ebre de Lie simple complexe et soit $\mathfrak{g} =
\mathfrak{n}_-\oplus
\mathfrak{h} \oplus \mathfrak{n}_+$, ``la'' d\'ecomposition triangulaire bien connue :
$\mathfrak{b} =
\mathfrak{h} \oplus \mathfrak{n}_+$ est une sous-alg\`ebre de Borel de $\mathfrak{g}$ et
$\mathfrak{n}_+$ en est le radical nilpotent. Panyushev, dans \cite{Pa}, examine cet
exemple
parmi d'autres (Remarks 1.5.1) et note que son in\'egalit\'e (avec $\mathfrak{a} =
\mathfrak{n}_+$,
id\'eal de $\mathfrak{b}$) s'\'ecrit : $\hbox{ind}(\mathfrak{b}) +
\hbox{ind}(\mathfrak{n}_+) \leq \dim(\mathfrak{b}/\mathfrak{n}_+) =
rg(\mathfrak{g})$, car
$\hbox{ind}(\mathfrak{b},\mathfrak{n}_+) = 0$ d'apr\`es un r\'esultat de A. Joseph. Il
dit aussi qu'il
peut \^etre ``conceptuellement'' prouv\'e qu'il y a \'egalit\'e : $\hbox{ind}(\mathfrak{b}) +
\hbox{ind}(\mathfrak{n}_+) = rg(\mathfrak{g})$.

        \medskip
        \item On reprend les notations du paragraphe 2. On va voir que l'\'egalit\'e de
Panyushev est
r\'ealis\'ee lorsqu'on prend $\mathfrak{g} = \mathfrak{n}$ et $\mathfrak{a} = z$ (qui n'est
commutatif que lorsque
$z$ est le centralisateur d'un \'el\'ement nilpotent principal). On a en effet :
$$
        \hbox{ind}\ \mathfrak{n} + \hbox{ind}\ z = 2\ \hbox{ind}\ \mathfrak{n} + \dim \delta
$$
$$
        \dim(\mathfrak{n}/z) + 2\ \hbox{ind}(\mathfrak{n} ,z) = \dim \delta + 2 \ \hbox{ind}\
\mathfrak{n}
$$
(On a appliqu\'e les formules de l'indice de Panyushev).

\smallskip
        En particulier, lorsque $z$ est commutatif, i.e. lorsque $z= \delta$, on a un cas
d'\'egalit\'e avec
$\mathfrak{a}$ commutatif. On constate que lorsqu'on prend $\mathfrak{a} = \delta$,
l'in\'egalit\'e
de Panyushev est stricte sauf lorsque $z=\delta$.

        \medskip
        \item On consid\`ere l'alg\`ebre de Lie $\mathfrak{g}$ du groupe de Mautner ;
$\mathfrak{g}$
admet $(P,Q,E,X)$ comme base avec :
$$
        [P;Q] = E,\quad [X,P] = Q, \quad [X,Q] = -P
$$
        \end{enumerate}

\bigskip
\noindent
$\centerdot$ L'indice de $\mathfrak{g}$ est \'evidemment 2. On consid\`ere l'id\'eal
$\mathfrak{a}$
admettant $(P,Q,E)$ comme base ; c'est l'alg\`ebre de Heisenberg de dimension 3 ; donc
$\hbox{ind}(\mathfrak{a}) = 1$ et ainsi :
$$
        \hbox{ind}(\mathfrak{g}) + \hbox{ind}(\mathfrak{a}) = 3.
$$

\bigskip
\noindent
$\centerdot$ Un calcul direct montre que $\hbox{ind}(\mathfrak{g},\mathfrak{a}) =
1$. Donc : 
$$
        \dim(\mathfrak{g}/\mathfrak{a}) + 2\ \hbox{ind}(\mathfrak{g},\mathfrak{a}) = 3.
$$

\vskip 7mm
\noindent
\textsc{Question 3} : Caract\'eriser les cas d'\'egalit\'e dans l'in\'egalit\'e de Panyushev.

        \vskip 12mm

 %%%%%%%%%%%%%%%%%%%%%%
\section{L'indice des contractions}
        
\noindent
\textbf{4.1.}~Soit $[\, , \,]$ une loi d'alg\`ebre de Lie sur un espace vectoriel $V$.
Pour chaque $g$
dans $GL(V)$ :
$$
        [x,y]_g = g([\bar{g}^1 .x,\ \bar{g}^1 y])\quad (x,y) \in V^2
$$
d\'efinit une (nouvelle) alg\`ebre de Lie sur $V$, et l'ensemble de ces structures
d'alg\`ebres de Lie est
l'orbite de $[\, , \, ]$ sous l'action (naturelle) de $GL(V)$ d\'ecrite plus haut.

\medskip
        Par d\'efinition, on appelle \underline{contraction} de l'alg\`ebre de Lie $[\, , \, ]$
toute alg\`ebre de
Lie qui appartient \`a l'adh\'erence de la $GL(V)$-orbite de $[\, , \, ]$. Pour ce qui
me concerne, il
s'agit de l'adh\'erence dans l'espace vectoriel complexe (ou r\'eel) constitu\'e par les
applications
bilin\'eaires altern\'ees de $V \times V$ dans $V$.

\medskip
        Il est imm\'ediat que l'indice ne varie pas le long de l'orbite de $[\, , \, ]$, et
que par contre, il se
peut qu'il augmente en certains points de l'adh\'erence de cette orbite. En tout cas :
$\hbox{ind}\,
\mathfrak{g} \leq \hbox{ind}\, \mathfrak{h}$ lorsque $\mathfrak{h}$ est une
contraction de
$\mathfrak{g}$.

\vskip 7mm
 %%%%
\noindent
\textbf{4.2.}~Soient $\mathfrak{g}$ une alg\`ebre de Lie et soit $\mathfrak{k}$ une
sous-alg\`ebre de
$\mathfrak{g}$, qui op\`ere naturellement dans $\mathfrak{g}/\mathfrak{k}$. Il est
bien connu que
l'alg\`ebre de Lie $(\mathfrak{g}/\mathfrak{k})\ltimes \mathfrak{k}$, produit
semi-direct relatif \`a
cette op\'eration de $\mathfrak{k}$ dans l'espace vectoriel $(\mathfrak{g}/\mathfrak{k})$
(consid\'er\'e comme alg\`ebre de Lie ab\'elienne), est une contraction de l'alg\`ebre de Lie
$\mathfrak{g}$ (dite d'Inonu-Wigner) (voir par exemple \cite{C}, proposition 1.3 (2)).

\vskip 7mm
 %%%%
\noindent
\textbf{4.3.}~Soit $\mathfrak{g} = \mathfrak{k}\oplus \mathfrak{p}$ une
d\'ecomposition de
Cartan d'une alg\`ebre de Lie \underline{simple r\'eelle}. L'alg\`ebre de Lie $\mathfrak{p}
\underset{\hbox{ad}}{\ltimes}\mathfrak{k}$, produit semi-direct de $\mathfrak{k}$ par
l'ab\'eliannis\'ee $\mathfrak{p}$, relativement \`a l'action adjointe de $\mathfrak{k}$ dans
$\mathfrak{p}$, est connue comme \'etant l'alg\`ebre de Lie du groupe des d\'eplacements
de Cartan.
D'apr\`es ce qu'on vient de dire, c'est une contraction de l'alg\`ebre de Lie
$\mathfrak{g}$.

\vskip 7mm
\noindent
\textbf{Proposition}~: \textit{Les alg\`ebres de Lie $\mathfrak{g}$ et $\mathfrak{p}
\underset{\hbox{ad}}{\ltimes}\mathfrak{k}$ ont m\^eme indice.}

\vskip 5mm
\noindent
\textsc{D\'emonstration} : Soit $\mathfrak{a}$ un sous-espace de Cartan de
$\mathfrak{p}$. On
note $G$ et $K$ les groupes associ\'es \`a $\mathfrak{g}$ et $\mathfrak{k}$
respectivement. On sait
que la dimension maximale des $K$-orbites dans $\mathfrak{p}$ (ou dans
$\mathfrak{p}^*$,
puisque la forme de Killing met $\mathfrak{p}$ en dualit\'e avec lui-m\^eme) est : $\dim
\mathfrak{p} - \dim
\mathfrak{a}$. On a donc :
$$
        \hbox{ind}(\mathfrak{k},\mathfrak{p}) = \dim \mathfrak{a}(= rg(G/K)).
$$
Soit $\mathfrak{m}$ l'annulateur dans $\mathfrak{k}$ d'un \'el\'ement r\'egulier de
$\mathfrak{a}$, et
soit
$\mathfrak{t}$ un tore maximal de $\mathfrak{m}$. On sait que : $\mathfrak{t} \oplus
\mathfrak{a}$ est une sous-alg\`ebre de Cartan de $\mathfrak{g}$, de sorte que :
$$
        rg(\mathfrak{g})(= \hbox{ind}\ \mathfrak{g}) = \dim \mathfrak{t} + \dim \mathfrak{a}
$$
avec : $\dim \mathfrak{t} = \hbox{ind}(\mathfrak{m})$ et $\dim \mathfrak{a} =
\hbox{ind}(\mathfrak{k},\mathfrak{p})$. Par ailleurs, d'apr\`es la formule de l'indice
des produits
semi-directs (\cite{Ra}) : 
$$
        \hbox{ind}(\mathfrak{p} \underset{\hbox{ad}}{\ltimes} \mathfrak{k}) =
\hbox{ind}(\mathfrak{m})
+ \hbox{ind}(\mathfrak{k},\mathfrak{p}).
$$
Donc : $\hbox{ind}(\mathfrak{g}) =
\hbox{ind}(\mathfrak{p}\underset{\hbox{ad}}{\ltimes}\mathfrak{k}$).

\vskip 7mm
 %%%%
\noindent
\textbf{4.4.}~\textsc{Remarques} : D'apr\`es une communication priv\'ee de Panyushev
(\cite{Pa 2}), ce r\'esultat est un cas particulier d'un fait plus g\'en\'eral depuis
longtemps connu de lui
:

\vskip 4mm
\noindent
        \textbf{Th\'eor\`eme} : \textit{Soient $\mathfrak{g}$ une alg\`ebre de Lie r\'eductive
complexe, $G$ le
groupe connexe associ\'e, $K$ un sous-groupe r\'eductif de $G$, $\mathfrak{k}$ l'alg\`ebre
de Lie de $K$,
et $\mathfrak{p}$ l'orthogonal de $\mathfrak{k}$ dans $\mathfrak{g}$ (relatif \`a une
forme bilin\'eaire
invariante). On a alors :}
$$
        \hbox{ind}(\mathfrak{p} \ltimes \mathfrak{k}) - \ \hbox{ind}\ \mathfrak{g} = 2 \
\hbox{comp.}(G/K)
$$
\textit{o\`u} $\hbox{comp.}(G/K)$ \textit{est la \underline{complexit\'e} de $G/K$.}

\medskip
\noindent
et qu'il dit n'avoir pas publi\'e explicitement. Par contre, une d\'emonstration adhoc
du corollaire
suivant se trouve dans (\cite{Pa}, theorem 3.5).

\vskip 4mm
\noindent
        \textbf{Corollaire} : \textit{Soit  $\mathfrak{g} = \mathfrak{k} \oplus
\mathfrak{p}$ une
$\mathbb{Z}_2$-graduation. Alors} $\hbox{ind}(\mathfrak{g}) =
\hbox{ind}(\mathfrak{p} \ltimes
\mathfrak{k})$.

\vskip 4mm
        Ainsi, lorsque $K$ est un sous-groupe r\'eductif de $G$, $\hbox{ind}(\mathfrak{g}) =
\hbox{ind}((\mathfrak{g}/\mathfrak{k})\ltimes \mathfrak{k})$ si et seulement si $K$
est un
sous-groupe sph\'erique de $G$.

\vskip 7mm
 %%%%
\noindent
\textbf{4.5.}~Voici un \underline{autre  exemple} de contraction ayant le m\^eme
indice que celui
de la contract\'ee. On reprend $\mathfrak{g} = \mathfrak{n}_- \oplus \mathfrak{h} \oplus
\mathfrak{n}_+$ comme dans 3.5 ci-dessus, et $\mathfrak{b} = \mathfrak{h}\oplus
\mathfrak{n}_+$. Alors la forme de Killing permet d'identifier
$(\mathfrak{g}/\mathfrak{b})^*$ \`a
$\mathfrak{n}_+$ et l'op\'eration naturelle de
$\mathfrak{b}$ dans $(\mathfrak{g}/\mathfrak{b})^*$ se lit comme \'etant l'op\'eration
adjointe de
l'alg\`ebre de Lie $\mathfrak{b}$ dans son id\'eal $\mathfrak{n}_+$. On sait que
$\hbox{ind}(\mathfrak{b},\mathfrak{n}_+) =0$, et pr\'ecis\'ement le stabilisateur
g\'en\'erique est
\underline{ab\'elien}. D'apr\`es la formule de l'indice (\cite{Ra}) :
$$
        \hbox{ind}((\mathfrak{g}/\mathfrak{b})\ltimes \mathfrak{b}) = rg(\mathfrak{g}) 
$$
puisque le stabilisateur g\'en\'erique est de dimension \'egale \`a : $\dim \mathfrak{b} - \dim
\mathfrak{n}_+ = rg(\mathfrak{g})$.

\vskip 7mm
\noindent
\textsc{Question 4} : Apr\`es un certain nombre de travaux consacr\'es \`a ``l'analogie''
entre l'analyse
harmonique sur $G$ et sur son contract\'e $\mathfrak{p} \ltimes K$ (repr\'esentations
unitaires
irr\'eductibles des 2 groupes, classes de conjugaison des 2 groupes$\ldots$) il serait
int\'eressant de
comprendre les relations entre les alg\`ebres d'invariants polynomiaux des 2 groupes,
et de m\^eme
que celles entre les centres des alg\`ebres enveloppantes respectives des 2 alg\`ebres
de Lie.

\vskip 7mm
\noindent
\textsc{Question 4 bis} : Est-ce que $\hbox{ind}((\mathfrak{g}/\mathfrak{p})\ltimes
\mathfrak{p}) 
= rg(\mathfrak{g})$ lorsque $\mathfrak{p}$ est une sous-alg\`ebre parabolique de
$\mathfrak{g}$ ?

\vskip 7mm
 %%%%
\noindent
\textbf{4.6.}~Il se peut que la remarque tout \`a fait \'el\'ementaire selon laquelle
l'indice d'une
contraction majore celui de la contract\'ee puisse \^etre utile pour calculer certains
indices.

\medskip
        Par exemple : Soit $z$ comme dans le paragraphe 2 ci-dessus ; on sait, d'apr\`es le
lemme de
Vinberg, que l'indice de $z$ est \underline{minor\'e} par le rang de $\mathfrak{g}$.
S'il existe une
contraction $\mathfrak{h}$ de $z$ dont l'indice est major\'e par le rang de
$\mathfrak{g}$, on aura
:
$$
        rg(\mathfrak{g}) \leq \hbox{ind}\ z \leq \hbox{ind}\ \mathfrak{h} \leq
rg(\mathfrak{g})
$$
d'o\`u $\hbox{ind}\ z = rg(\mathfrak{g})$, ce qui assure dans ce cas la validit\'e de la
conjecture dite
d'Elashvili. Cela suppose \'evidemment de conna”tre \underline{des} contractions de
$z$ et de savoir
calculer leurs indices. Il y a au moins une possibilit\'e, chaque fois qu'on dispose
d'une sous-alg\`ebre
de Lie $\mathfrak{k}$ de $z$ : on prend $\mathfrak{h} = (z/\mathfrak{k})\ltimes
\mathfrak{k}$ et
on peut utiliser la formule de
\cite{Ra} pour en calculer l'indice. Plus particuli\`erement, d'apr\`es un th\'eor\`eme
classique de
Kostant, il existe une sous-alg\`ebre ab\'elienne $\mathfrak{k}$ de dimension $r$ (o\`u on
note $r$ le rang
de $\mathfrak{g}$), et dans ce cas, l'application de la formule de l'indice s'en
trouve facilit\'ee (le
stabilisateur g\'en\'erique est ab\'elien, de sorte que son indice co\"\i ncide avec sa
dimension).

\vskip 5mm
\noindent
\textsc{Exemple} : On suppose $e$ ``sous-r\'egulier'' de sorte que $\dim z = r+2$. On
prend une
sous-alg\`ebre ab\'elienne $\mathfrak{k}$ de $z$, de dimension $r$, de sorte que
$z/\mathfrak{k}$ est
de dimension 2. Les orbites de dimension maximale $d$ pour l'op\'eration naturelle de
$\mathfrak{k}$ dans le dual de
$(z/\mathfrak{k})$ sont \'evidemment telles que : $d = 0,1$ ou 2.

\vskip 5mm
        $\centerdot$ Lorsque $d=2$, l'indice de la repr\'esentation de $\mathfrak{k}$ dans
$(z/\mathfrak{k})$ est z\'ero, tandis que le stabi\-lisateur g\'en\'erique dans
$\mathfrak{k}$ est de
dimension $(r-2)$. Donc : $\hbox{ind}((z/\mathfrak{k}\ltimes \mathfrak{k}) = r-2$. 
Ce cas est impossible.

\vskip 5mm
        $\centerdot$ Lorsque $d=1$, le m\^eme raisonnement donne :
$$
        \hbox{ind}((z/\mathfrak{k})\ltimes \mathfrak{k}) = r.
$$

\vskip 5mm
        $\centerdot$ Lorsque $d=0$, le m\^eme raisonnement donne :
$$
        \hbox{ind}((z/\mathfrak{k})\ltimes \mathfrak{k}) = r+2.
$$

        Le troisi\`eme cas est celui o\`u $\mathfrak{k}$ est un id\'eal de $z$. En
\underline{conclusion} : si
$\mathfrak{k}$ est une sous-alg\`ebre ab\'elienne de $z$, de dimension $r$, et qui n'est
pas un id\'eal de
$z$, alors :
$\hbox{ind}((z/\mathfrak{k})\ltimes \mathfrak{k}) = r$.

\medskip
        Le fait que l'indice de $z$, dans ce cas ($e$ sous-r\'egulier), vaut $r$ a \'et\'e
d\'emontr\'e en 2 lignes par
Panyushev (3.4, Corollary dans \cite{Pa}). Noter toutefois qu'il utilise la
caract\'erisation des
nilpotents principaux par la commutativit\'e de leurs centralisateurs.

\vskip 7mm
\noindent
\textsc{Question 5} : Existe-t-il toujours une sous-alg\`ebre ab\'elienne, de dimension
$r$, qui ne soit
pas un id\'eal de $z(e)$, $e$ \'etant un nilpotent sous-r\'egulier ?

\vskip 12mm
 %%%%%%%%%%%%%%%%%%%
\section{Des exemples d'additivit\'e de l'indice}

\bigskip
        Soit $\mathfrak{g} = \mathfrak{g}_0 \oplus \mathfrak{g}_1$ une d\'ecomposition d'une
alg\`ebre de
Lie en somme directe de 2 sous-alg\`ebres de Lie $\mathfrak{g}_0$ et $\mathfrak{g}_1$.
On donne
dans ce paragraphe des exemples de cette situation o\`u $\hbox{ind}(\mathfrak{g}) =
\hbox{ind}(\mathfrak{g}_0) + \hbox{ind}(\mathfrak{g}_1)$.

\vskip 7mm
%%%%
\noindent
\textbf{5.1.}~Soient $e$ un \'el\'ement nilpotent dans une alg\`ebre de Lie simple
$\mathfrak{g}$ et
$z$ le centralisateur de $e$. Supposons que $z$ soit Lie-compl\'ement\'e, i.e. qu'il
existe une
sous-alg\`ebre $\mathfrak{g}_1$ de $\mathfrak{g}$ telle que : $\mathfrak{g} = z \oplus
\mathfrak{g}_1$ (des exemples de cette situation ont \'et\'e d\'ecrits dans \cite{Ra 2} et
dans
\cite{Sa}). 

        On sait, dans ce cas, que $\mathfrak{g}_1$ est d'indice z\'ero. Par cons\'equent
``$\hbox{ind}(\mathfrak{g}) = \hbox{ind}\ z + \hbox{ind}\ \mathfrak{g}_1$'' \'equivaut
 \`a :
``$\hbox{ind}\ z = rg(\mathfrak{g})$'' i.e., \'equivaut \`a la validit\'e de la conjecture
d'Elashvili.
Lorsque $e$ est un \'el\'ement d'une orbite sph\'erique, on a bien : $\hbox{ind}(z) =
rg(\mathfrak{g})$
(\cite{Pa}, theorem 3.5) et (\cite{Pa 3}, theorem 3.5), et d'apr\`es (\cite{Sa},
th\'eor\`eme 3.3) $z$ est
Lie-compl\'ement\'e de sorte que :

\bigskip
        Lorsque $e$ appartient \`a une orbite sph\'erique, il existe des sous-alg\`ebres
$\mathfrak{g}_1$ telles
que : $\mathfrak{g} = z \oplus \mathfrak{g}_1$ et pour chaque telle sous-alg\`ebre
$\mathfrak{g}_1$, on a :
$$
        \hbox{ind}\ \mathfrak{g} = \hbox{ind}\ z + \hbox{ind}\ \mathfrak{g}_1.
$$

\vskip 7mm
 %%%%
\noindent
\textbf{5.2.}~Pour des raisons de commodit\'e, $\mathfrak{g} = \mathfrak{g} \ell(n)$
(le cas $sl(n)$
s'en d\'eduit). Soit $p$ un diviseur de $n$, et soit $e = J^p$ o\`u :
$$
        J = \begin{pmatrix}
0        &0                &                &        &0\\
1        &0        &        &        &\vdots\\
0        &1        &\ddots        &                &\vdots\\
\vdots                &        &\ddots        &\ddots\\
0        &0        &\hdots        &1        &0\\
\end{pmatrix}
$$
est la matrice de Jordan de taille $n$. D'apr\`es \cite{Ra 3}, l'alg\`ebre de Lie
$\mathfrak{g}_1$,
constitu\'ee par les matrices ayant leurs $p$ derni\`eres lignes nulles, est d'indice
z\'ero, et le
centralisateur $z$ de $e$ dans $\mathfrak{g}\ell(n)$ est l'alg\`ebre de Lie constitu\'ee
par les
matrices dont l'\'ecriture en blocs $p \times p$, est la suivante :
 $$
\begin{pmatrix}
A_1                                        &                                                         &                                                                                                &0                                                        &\\
                                                                                        &A_1                        &                                                                                        &                                                                        &\\
A_2                                        &                                                                        &\ddots        &                                                                &                \\
\vdots                &A_2                &                                                                                                &\ddots                &\\
\vdots                &                                                                &\ddots                &                                                                &\ddots\\
A_k                                        &\ldots        &                                                                                        &A_2                        &                        &A_1\\
\end{pmatrix}\qquad (n = kp)
$$
o\`u les $A_i$ sont des matrices $p\times p$ quelconques. On a alors :
$\mathfrak{g}\ell(n) = z
\oplus \mathfrak{g}_1$. On voit facilement que l'application qui associe \`a chaque
matrice, du type
pr\'ec\'edent, l'expression : $A_1 + tA_2 +\cdots + t^{k-1}A_k$ est un isomorphisme de
l'alg\`ebre de
Lie $z$ sur une alg\`ebre de Takiff g\'en\'eralis\'ee, du type de celles \'etudi\'ees dans
\cite{R-T}. Du
th\'eor\`eme 2.8 de l'article cit\'e, on d\'eduit : $\hbox{ind}(z) = kp = n$. Ainsi
$\hbox{ind}(z) = rg(\mathfrak{g})$ pour ce type de centralisateur.

\medskip
        On notera que d'apr\`es (\cite{Pa}, 3.6), la conjecture $\hbox{ind}(z) =
rg(\mathfrak{g})$ aurait \'et\'e
d\'emontr\'ee dans beaucoup de cas (y compris celui de $sl(n)$) par Elashvili, mais que
ses calculs
n'ont pas \'et\'e publi\'es.

\vskip 7mm
 %%%%
\noindent
\textbf{5.3.}~Ici $\mathfrak{g} = \mathfrak{g}\ell(n),\ \mathfrak{g}_0 = so(n),\
\mathfrak{g}_1
= \mathfrak{b}$ (l'alg\`ebre de Lie des matrices triangulaires sup\'erieures), de sorte que
$\mathfrak{g} = \mathfrak{g}_0 \oplus \mathfrak{g}_1$. On sait que :

$$
        \hbox{ind} \ \mathfrak{g}_0 = \Big[{n \over 2}\Big] \ \hbox{et}\ \hbox{ind} \
\mathfrak{g}_1 =
\Big[{{n+1} \over 2}\Big]
$$
et il vient imm\'ediatement : $\hbox{ind}\ \mathfrak{g}_0 + \hbox{ind}\ \mathfrak{g}_1 =
n = rg(\mathfrak{g})$.

\vskip 7mm
\noindent
\textsc{Question 6} : Soit $\mathfrak{g} = \mathfrak{k}\oplus \mathfrak{p}$ une
d\'ecomposition de
Cartan comme dans 4.3 ci-dessus. On sait que $\mathfrak{g} = \mathfrak{k} \oplus
\mathfrak{b}$
avec
$\mathfrak{b} = \mathfrak{a} \oplus \mathfrak{n}$ (les notations sont classiques),
i.e. :
$\mathfrak{g} = \mathfrak{k} \oplus \mathfrak{a} \oplus \mathfrak{n}$ est une
d\'ecomposition 
d'Iwasawa de $\mathfrak{g}$.

\vskip 4mm
        A-t-on  : $\hbox{ind}\ \mathfrak{g} = \hbox{ind}(\mathfrak{k}) +
\hbox{ind}(\mathfrak{b})$ ?

\vskip 7mm
 %%%%
\noindent
\textbf{5.4.}~\begin{enumerate}
\item Soit $\mathfrak{g} = \mathfrak{n}_-\oplus \mathfrak{h}\oplus \mathfrak{n}_+$ une
d\'ecomposition triangulaire de
$\mathfrak{g}$, comme dans 3.5 ci-dessus. On a alors $\mathfrak{g} =
\mathfrak{n}_-\oplus
\mathfrak{b}$, avec $\mathfrak{b} = \mathfrak{h}\oplus \mathfrak{n}_+$, et comme
d\'ej\`a dit dans
3.5.a, (\cite{Pa}, 1.5.1), on a : $\hbox{ind}\ \mathfrak{g} =
\hbox{ind}(\mathfrak{b}) +
\hbox{ind}(\mathfrak{n}_+)$, donc : 
$$
        \hbox{ind}(\mathfrak{g}) = \hbox{ind}(\mathfrak{n}_-) + \hbox{ind}(\mathfrak{b}).
$$

\medskip
\item On notera que  cet exemple est un cas particulier de la situation \'etudi\'ee par
B. Kostant
(\cite{Ko}) : soit $\mathfrak{g}$ une alg\`ebre de Lie semi-simple (sur $\mathbb{R}$ ou
$\mathbb{C}$), et soit $\Phi$ une forme bilin\'eaire non d\'eg\'en\'er\'ee, sym\'etrique et
invariante sur
$\mathfrak{g}$ (disons la forme de Killing de $\mathfrak{g}$). On note $\theta$ une
involution de
Cartan de
$\mathfrak{g}$, ce qui permet de d\'efinir sur $\mathfrak{g}$ un produit scalaire :
$$
        <x|y> \, = \Phi(x,-\theta(y))\qquad (x,y) \in \mathfrak{g} \times \mathfrak{g}.
$$
Kostant dit d'une sous-alg\`ebre de Lie $\mathfrak{a}$ de $\mathfrak{g}$ qu'elle est
``Lie-sommante'' lorsque son orthogonal (pour la forme de Killing) $\mathfrak{a}^o$
est une
sous-alg\`ebre de Lie de $\mathfrak{g}$ ou, ce qui revient au m\^eme, lorsque son
orthogonal (pour le
produit scalaire introduit) $\mathfrak{a}^\perp$ est lui-m\`eme une sous-alg\`ebre de
Lie de
$\mathfrak{g}$, de sorte qu'on a alors : $\mathfrak{g} = \mathfrak{a} \oplus
\mathfrak{a}^\perp$. Kostant remarque que lorsque $\mathfrak{a}$ est une sous-alg\`ebre
parabolique $\mathfrak{p}$ de $\mathfrak{g}$, son orthogonal $\mathfrak{a}^o$ est
son radical
nilpotent $\mathfrak{p}^u$ (et $\mathfrak{a}^\perp = \theta(\mathfrak{p}^u)$), de
sorte que :
$\mathfrak{g} = \mathfrak{p} \oplus \theta(\mathfrak{p}^u)$ pour toute sous-alg\`ebre
parabolique $\mathfrak{p}$ de $\mathfrak{g}$. Le cas o\`u $\mathfrak{p} =
\mathfrak{b}$ est une
sous-alg\`ebre de Borel de $\mathfrak{g}$ a \'et\'e pr\'esent\'e ci-dessus, et on a bien :
$\hbox{ind}(\mathfrak{g}) = \hbox{ind}(\mathfrak{b}) + \hbox{ind}(\mathfrak{b}^u)$.

\medskip
\item \underline{Un autre exemple de cette situation} : ici $\mathfrak{g} = sl(n)$ et
$\mathfrak{p}$ est l'alg\`ebre de Lie constitu\'ee par les matrices $(x_{ij})$ de trace
nulle et telles que
$x_{n,1} = x_{n,2} =\cdots = x_{n,n-1} = 0$ ; l'alg\`ebre $\mathfrak{p}$ est une
sous-alg\`ebre
parabolique maximale de $\mathfrak{g}$, et la d\'ecomposition ``orthogonale'' de
Kostant qui lui
est associ\'ee est : $sl(n) = \mathfrak{p} \oplus \mathfrak{q}$, o\`u $\mathfrak{q}$ est
l'ensemble
des matrices $y = (y_{ij})$ telles que : $y_{ij} = 0$ lorsque $1 \leq i \leq n-1$ et
$y_{nn} = 0$. On
remarquera que dans ce cas, l'alg\`ebre de Lie $\mathfrak{q}$ est ab\'elienne, de sorte
que :
$\hbox{ind}(\mathfrak{q}) = \dim \mathfrak{q} = n-1$. Comme $\mathfrak{p}$ est
d'indice z\'ero
(c'est un cas particulier de ce qui est contenu dans 5.2 ci-dessus), on a ici aussi :
$\hbox{ind}(\mathfrak{g}) = \hbox{ind}\ \mathfrak{p} + \hbox{ind}(\mathfrak{p}^u)$.
        \end{enumerate}

\vskip 7mm
\noindent
\textsc{Question 7} : Pour quel type de parabolique a-t-on :
$\hbox{ind}(\mathfrak{g}) =
\hbox{ind}\ \mathfrak{p} + \hbox{ind}\ \mathfrak{p}^u$ ?

\vskip 7mm
 %%%%
\noindent
\textbf{5.5.}~Il n'est pas vrai en g\'en\'eral que lorsque $\mathfrak{g} =
\mathfrak{g}_0 \oplus
\mathfrak{g}_1$ est une d\'ecomposition en somme de 2 sous-alg\`ebres de Lie, m\^eme si plus
particuli\`erement il s'agit d'une d\'ecomposition \`a la Kostant, on ait : $\hbox{ind}\
\mathfrak{g} =
\hbox{ind} \ \mathfrak{g}_0 + \hbox{ind} \ \mathfrak{g}_1$. Prenons $\mathfrak{g} =
\mathfrak{g}\ell(4)$ pour $\mathfrak{g}_1$ l'alg\`ebre de Lie des matrices dont les 2
derni\`eres
lignes sont nulles, et pour $\mathfrak{g}_2$ celle des matrices dont les 2 premi\`eres
lignes sont
nulles. On a : $\mathfrak{g} = \mathfrak{g}_1 \oplus \mathfrak{g}_2$ et
$\hbox{ind}(\mathfrak{g}_1) = \hbox{ind}(\mathfrak{g}_2) = 0$ (\cite{Ra} 2.16).

\vskip 12mm

\renewcommand{\refname}{Bibliographie} %%%pour m\'emoire : \cite{ } dans le texte

\end {document}